\documentclass[notitlepage, 11pt]{article}
\usepackage{amsmath}
\usepackage{amsfonts}
\usepackage{amssymb}
\usepackage{latexsym}

\usepackage{color}

\newtheorem{lemma}{Lemma}
\newtheorem{theorem}{Theorem}

\newtheorem{definition}{Definition}

\newtheorem{condition}{Condition}

\newtheorem{remark}{Remark}

\newcommand{\Z}{\mathcal{Z}}
\newcommand{\M}{\mathcal{M}}
\newcommand{\N}{\mathbb{N}}
\newcommand{\Li}{\mathcal{L}}
\newcommand{\Sp}{\mathcal{S}}
\newcommand{\eps}{\epsilon}
\newcommand{\noi}{\noindent}
\newcommand{\spa}{\vspace{.2in}}

\newcommand{\R}{\mathbb{R}}

\newcommand{\p}{\mathbb{P}}

\newcommand{\la}{\lambda}

\newcommand{\Ups}{\Upsilon}
\newcommand{\del}{\delta}
\newcommand{\Om}{\Omega}

\newcommand{\be}{\begin{equation}}
\newcommand{\ee}{\end{equation}}
\newcommand{\B}{\mathcal{B}((0,\infty))}
\newcommand{\BB}{\mathcal{B}([0,\infty))}
\newcommand{\lan}{\langle}
\newcommand{\ran}{\rangle}

\title{Fluid limits of many-server queues with state\\
dependent service rates}
\author{Anup Biswas\\ \\
Department of Electrical Engineering\\
Technion--Israel Institute of Technology\\
Haifa 32000, Israel}

\date{}

\begin{document}
\maketitle

\begin{abstract}
We study a many-server queuing system with general service time distribution and state dependent service rates. The dynamics of the system are
modeled using measure valued processes which keep track of the residual service times. Under suitable conditions, we prove
the existence of a unique fluid limit.

\spa

\noi{\bf AMS subject classifications:}\,\, 60F17, 60K25, 90B22

\spa

\noi{\bf Keywords:}\,\,
Many-server queues, GI/G/n queue, fluid limits, measure-valued processes, state dependent service rate, call centers.
\end{abstract}

\section{Introduction}

In recent days, many-server queuing systems have received much attention due to its applications to call center. Thus it has become important to study its asymptotic properties
to gain insight into the behavior of these systems. Studying different scaling limits (fluid or diffusion scaling) are established tradition in queuing theory. In the celebrated
work of Halfin and Whitt \cite{halfin-whitt}, it was shown that  for Poissonian arrival and exponential service time, with positive probability, there is a positive queue in the asymptotic regime if the arrival rate $\la_n$ and the number of servers $n$ both goes to infinity in a manner that $\la_n=n-\beta\sqrt n$ for some $\beta>0$. Fluid and diffusion limits for the
total number of customers in a network with time varying Poissonian arrival and stuffing was obtained in \cite{mandel-mass-reim}. This work was generalized in \cite{liu-whitt} for
$G_t/M_t/s_t$ queuing networks with abandonment where the authors studied long-time behavior of the fluid limits.

\spa

A recent statistical study by Brown et. al. \cite{brown} suggests that it may be more appropriate to consider non-exponential service times. In particular, it is log-normal in some
cases as shown in \cite{brown}. This emphasizes the need to consider many-server model with generally distributed service times. In \cite{whitt}, Whitt considered G/G/n network with
abandonment and proposed a deterministic fluid approximation. The author proved convergence for discrete time model. In \cite{kaspi-ramanan}, Kaspi and Ramanan considered G/G/n model
and obtained a measure-space valued fluid limit. Later Kang and Ramanan generalized this work by allowing the customers abandonment in \cite{kang-ramanan}. In \cite{atar-kaspi-shim},
Atar et al. studied multi-class many-server queues with fixed priority and established the existence of unique fluid limits. Kang and Ramanan studied ergodic properties of the scaled
processes for GI/G/n+G model and its relation with invariant states of the fluid limit in \cite{kang-ramanan1}. Reed in \cite{reed}, established the fluid and diffusion limits
of the customer-count processes for many-server queuing system under the finite first moment assumption on service time distribution. In \cite{zhang}, Zhang obtained the fluid limits for GI/G/n+G queuing systems.

\spa

All of the above models consider servers that serve the customers at a constant rate $1$. In this work, we allow the servers to adjust their service rate depending on the number
of customers in the system (equivalently, the number of customers in the queue). It is often useful to increase the service rates when the queue length is large.
Management may also be interested to adjust the service rate depending on customers feedback. State dependent arrival and service rate were first introduced in \cite{yamada} for
conventional heavy traffic approximations. In case of single server models, processor sharing model is an example where service rate at any instant of time depends on the number of customers in the system. For some recent developments on processor sharing, we refer \cite{gro-rob-zwart}, \cite{zhang-dai-zwart}, \cite{ramanan-reiman}.

\spa

In this work, we consider a system with $n$-homogeneous servers.
Customers arrival is given by a renewal process and
customers are served under FCFS policy. Arrived customers do not leave the system until served.
Let $X^n_t$ denote the number of customers in the $n$-th system at time $t\geq 0$. Define $\bar X^n_t=\frac{1}{n}X^n_t$.
The service rate of each server at time $t$ is given by $k^n(\bar X^n_t)$ for some bounded map $k^n$ on $[0,\infty)$. Note that, $k^n(\bar X^n_t)$ could be $0$ for non-zero
 $\bar X^n_t$. The system is described by $(Q^n_t,\Z^n_t)$ where $Q^n_t$ denotes the number of customers waiting in the queue at time $t$ and $\Z^n_t$ is a non-negative Borel measure
 on $(0,\infty)$ such that $\Z^n_t(C)$ denotes the number of customers in service with their remaining service
 requirements in $C$, for $C\in\B$.
 Thus $\Z^n$ is a measure-valued process that keeps track of customers remaining service requirements. Measure-valued processes
that keep track of residual service requirements of individual customers, have been considered earlier in literature (see \cite{dec-pascal}, \cite{gro-rob-zwart}, \cite{zhang-dai-zwart},
\cite{zhang}).
In \cite{kaspi-ramanan}, the authors used measure-valued processes that keep track of the time spent by the individual customers in service. Also their proof relies on the fact that
there exists a compensator for the departure processes (see Corollary 5.5 there). Since we are allowing our service rate to be dynamic depending on $\bar X^n$, getting an explicit compensator for analogous processes as in \cite{kaspi-ramanan}
is a hard problem. In this work, we closely follow the approach in \cite{zhang-dai-zwart}
(see also \cite{zhang}).
We show that the fluid limit of $\frac{1}{n}\Z^n_t$ is uniquely determined by an integral equation, referred as fluid model equation(\eqref{theo2} below). A similar type of equation was also obtained in \cite{zhang-dai-zwart}.

\spa

This paper is organized as follows. In Section 2, we introduce our $GI/G/n$ model with our basic assumptions and state our main result. The uniqueness part
of the main result is done in Section \ref{fluid}. In Section 3, we prove relative compactness of the stochastic pre-limits and characterizations of the limits are done in Section 4.
Finally, in Appendix we prove existence result for solution to certain integral equation and recall some results from \cite{zhang-dai-zwart} those are used in this paper.

\subsection{Notations}
The following notations will be used throughout this paper. By $\N, \R,$ we denote the set of natural numbers and the
set of real numbers, respectively. Given $a,b\in\R$, the maximum (minimum) is denoted by $a\vee b$($a\wedge b$). We
use $a^+$ for $a\vee 0$. We define $\R_+=[0, \infty)$.
For any $A\subset [0,\infty)$, we define $A^\eps=\{x\geq 0: \inf_{a\in A}|x-a|<\eps\}$. For any $x\in\R$, the sets $(x, \infty), [x, \infty)$
will be denoted by $C_x, \bar C_x$ respectively.
For any topological space $\Sp$,
$C_b(\Sp)$ denotes the set of all real valued bounded, continuous map on $\Sp$ and $\mathcal{B}(\Sp)$ is used to denote the Borel $\sigma$-field of $\Sp$. $C([a,b], \Sp)$ will denote the set of all continuous function from $[a,b]$ to $\Sp$. For any $f\in C([a,b], \R)$, we define $|f|_{ab}=\sup_{x\in[a,b]}|f(x)|$.

The set of all non-negative finite Borel measures on $[0,\infty)$ is denoted by $\M$ and $\M_+$ denotes the subset of $\M$
containing all the measures having no atom at $\{0\}$. For any $\mu\in\M$ and Borel measurable function $g$ on $[0,\infty)$, we define $\langle g, \mu\rangle=\int g d\mu$.
For $\mu_1,\mu_2\in\M$, the Prohorov metric is defined by
$$\rho(\mu_1,\mu_2)=\inf\{\eps>0 : \mu_1(A)\leq \mu_2(A^\eps)+\eps, \mu_2(A)\leq \mu_1(A^\eps)+\eps,
\mbox{for all closed}\ A\subset[0,\infty)\}.$$
It is well known that $(\M, \rho)$ is a Polish space (see Appendix in \cite{delay-verejones}). Also this topology is equivalent
to the weak topology on $\M$ which is characterized as follows: $\mu_n\rightarrow \mu$ in weak topology if and only
if
\begin{center}
$\langle f, \mu_n\rangle\to\langle f, \mu\rangle$ for all $f\in C_b([0,\infty)$.
\end{center}

\noi Given any Polish space $(E, \pi)$, $D([0,\infty), E)$ denote the
space of functions that are right-continuous with finite left limits (RCLL). Endow the space $D([0,\infty), E)$ with the
Skorohod-Prohorov-Lindvall metric or $J_1$ metric,
defined as
$$
 d(\phi,\phi^\prime)=\inf_{f\in\Ups}\Big(\|f\|^\circ\vee\int_0^\infty e^{-u} d_u(\phi,\phi^\prime, f)du \Big),
\qquad \phi, \phi^\prime\in D([0, \infty), E)
$$
where
$$
d_u(\phi,\phi^\prime, f)=\sup_{t\geq 0}[\pi(\phi(t\wedge u),\phi^\prime(f(t)\wedge u))\wedge 1],
$$
and $\Ups$ is the set of strictly increasing, Lipschitz continuous functions from $[0, \infty)$ onto itself,
with
$$\| f\|^\circ=\sup_{0\leq s<t}\Big|\log\frac{f(t)-f(s)}{t-s}\Big|<\infty.$$
As is well known \cite{ethier-kurtz}, $D([0, \infty), E)$ is a Polish space under the induced topology.

We use $"\Rightarrow"$ to denote the convergence in the sense of distribution.

\section{Queuing model}
In this section, we describe our GI/G/n model and the measure valued state descriptors.
We assume that for each $n$, all the stochastic variables below, are defined on probability space $(\Omega_n, \mathcal{F}_n,\p_n)$.
The system contains
$n$ identical servers. Each arriving customer has a single service requirement and is served by a single server.
The customers are served by FCFS policy and they leave the system once their service is completed. We do not allow
the customers to renege the system until their job is done. We also assume that the system works under work conserving policy i.e., all the servers are busy if there is a queue. We assume the following:
\begin{itemize}
\item Customers arrive according to a renewal process $E^n_t$ with mean inter-arrival time $\frac{1}{\la^n}$ for some
$\la^n>0$.
\item Service requirement of the arriving $i$-th customer is given by $v^n_i$ where $\{v^n_i\}_{i=-\infty}^{i=\infty}$ is an positive valued
i.i.d. sequence with common distribution $\nu^n$.
\end{itemize}

Upon arrival the customers join the queue if all the servers are busy.
At time $t,$ all the servers serve the customers at a rate $k^n(\bar{X}^n_t)$ where $\bar X^n_t=\frac{X^n_t}{n}$ and
$X^n_t$ denotes the number of customers in system at time $t$. It is reasonable to assume that $k^n:\R_+\to\R_+$ is a bounded, Borel measurable function.
 Let $\tau^n_i, i\in\N,$ be the time when the $i$-th customer starts its service.
Then for $t\geq 0$, the remaining service of $i$-th customer is given $v^n_i-\int^t_{\tau^n_i}k^n(\bar{X}^n_s)ds$
(non positive quantity implies that customer's job is completed) provided $\tau^n_i\leq t$. We use negative indices to denote the customers
in system at time $t=0$.
Let $X^n_0$ denote the number customers at
time $t=0$ with remaining job $\tilde v^n_i, i= -X^n_0+1,\ldots, -[X^n_0-n]^+,$ for the customers in service at time $t=0$ where
$\tilde v_i^n, i=-X^n_0+1,\ldots, -[X^n_0-n]^+,$ are some random variables defined on $(\Om_n,\mathcal{F}_n,\p_n)$.
Also let $Q^n_0$ denote the number of customers in the queue at time $t=0$. Hence $Q^n_0=[X^n_0-n]^+$.

\spa

For $t\geq 0$, $\Z^n_t$ denotes a measure in $\M_+$ such that $\Z^n_t(C)$ denotes the number of customers in service with remaining service
requirement in $C$ for $C\subset ((0,\infty))$. Hence the total number of customers in service at time $t$ is given by $Z^n_t=\Z^n_t(0,\infty)$. Let
$Q^n_t$ be the number of customers waiting in the queue. Define $B^n_t=E^n_t-Q^n_t$. Then $B^n_t+1$ denotes the index of the head of
the customers in the queue waiting to be served. For $0\leq s\leq t\leq T$, we define $S^n(s, t)=\int_s^tk^n(\bar X^n_u)du$.
Hence a precise description of $\Z^n$ is given by

\be\label{1}
\Z^n_t(C)=\sum_{i=-X^n_0+1}^{-Q^n_0}\del_{\tilde v^n_i}(C+S^n(0, t))+ \sum_{i=-Q^n_0+1}^{B^n_t}\del_{v^n_i}(C+S^n(\tau^n_i, t)),
\ee
for $C\in\B$. Additional obvious relation satisfied by the processes, for $t\geq 0$, are as follows:
\be\label{2}
X^n_t = Q^n_t+Z^n_t,
\ee
\be
Q^n_t = [X^n_t-n]^+. \label{3}
\ee
It is easy to see from \eqref{2} and \eqref{3} that $Z^n_t=X^n_t\wedge n$. We extend $\Z^n_t$ on $\BB$ by setting $\Z^n_t(\{0\})=0$.
It is easy to see that $\Z^n$ takes values in $D([0, \infty),\M)$.

\subsection{The fluid model}\label{fluid}
In this section, we define the fluid model and state our main theorem. We also state the set of assumptions that are
used to prove this result. We define
$$ \bar X^n_t=\frac{X^n_t}{n}, \bar\Z^n_t=\frac{1}{n}\Z^n_t, \bar Q_t^n=\frac{Q^n_t}{n}, \bar Z_t^n=\frac{Z^n_t}{n}.$$
The fluid scaling of arrival process $E^n_t$ is define as $\bar E^n_t=\frac{1}{n}E^n_t.$ Thus the fluid pre-limit equations
are (in analogy with \eqref{1}--\eqref{3}) given by
\be\label{4}
\bar\Z^n_t(C)=\frac{1}{n}\sum_{i=-X^n_0+1}^{-Q^n_0}\del_{\tilde v^n_i}(C+S^n(0, t))+ \frac{1}{n}\sum_{i=-Q^n_0+1}^{B^n_t}\del_{v^n_i}(C+S^n(\tau^n_i, t)),
\ee
for $C\in\BB$.
\be\label{5}
\bar X^n_t = \bar Q^n_t+ \bar Z^n_t,
\ee
\be
\bar Q^n_t = [\bar X^n_t-1]^+. \label{6}
\ee

We assume the following conditions:
\begin{condition}\label{cond-1}
\begin{itemize}
\item[(a)] $\frac{\la^n}{n}\to\la$ for some $\la\in[0,\infty)$ and $\bar E^n_t\Rightarrow \lambda t$\ in the sense of distribution in $D([0, \infty), \R_+)$;
\item[(b)] There exists a probability measure $\nu$ with bounded, Lipschitz continuous, density $g:\R_+\to\R_+$
 such that $\nu^n\rightarrow\nu$ as $n\to\infty;$
\item[(c)] $(\bar\Z^n_0, \bar Q^n_0)\Rightarrow (\Z_0, Q_0)$ in $\M\times\R_+$ as $n\to\infty$ where the function $F(x)=\Z_0([x,\infty))$ is Lipschitz continuous
on $\R_+$ and $(\Z_0, Q_0)$ is a deterministic element in $\M\times\R_+$.
\end{itemize}
\end{condition}
Since $\la t$ is a continuous, deterministic path, one obtains the convergence of the scaled arrival process in probability i.e., for any $T, \eps>0$
\be\label{7}
\lim_{n\to\infty}\p_n(\sup_{0\leq t\leq T}|\bar E^n_t-\la t|>\eps)=0.
\ee
We impose the following condition on the state dependent service rates:
\begin{condition}\label{cond-2}
There exists a bounded, Lipschitz continuous map $k:\R_+\to\R_+$ such that $k^n\to k$ as $n\to\infty$ uniformly on compact subsets of $\R_+$.
\end{condition}
Let $G(\cdot)$ be the distribution function of $\nu$ i.e., $G(x)=\nu([0,x])$. Define $G^c(x)=1-G(x)$. Also from \eqref{5} , \eqref{6}
and Condition \ref{cond-1}(c) the followings hold:
$$ \lan 1,\Z_0\ran=Z_0,\  X_0 =  Q_0+  Z_0,\  Q_0 = [ X_0-1]^+. $$

\begin{theorem}\label{main}
Assume Conditions \ref{cond-1}, \ref{cond-2} to hold. Then as $n\to\infty$, $(\bar\Z^n,\bar Q^n)\Rightarrow (\Z, Q)$ in $D([0,\infty), \M\times\R_+)$
where $(\Z, Q)$ is uniquely defined by the followings:
\begin{itemize}
\item for all $t\geq 0$,
\be\label{theo1}
 \Z_t((0,\infty))=Z_t,\  X_t =  Q_t+  Z_t,\  Q_t = [ X_t-1]^+,
\ee
\item for all $t\geq 0$ and $x\geq 0$,
\be\label{theo2}
\Z_t(\bar C_x)=F(x+S(0,t))+\int_0^tG^c(x+S(s,t))dB_s,
\ee
where $S(s,t)=\int_s^tk(X_u)du$ and $B_t=\la t-Q_t$. We refer \eqref{theo2} as fluid model equation and will be denoted by $(k, \la, \nu)$.
\end{itemize}

\end{theorem}

\noi{\bf Proof:}\ The existence of a limit satisfying the above properties will be proved in Section 3 and 4. So it is enough to prove the uniqueness of the limit here.
Since $B^n_t$ is nondecreasing and we can have pointwise convergence for each subsequential limit(by Skorohod representation theorem), $B_t$ is also non-decreasing. Therefore
\eqref{theo2} makes sense. Also for any $a>0$,
$$\Z_t([0,a))\leq \mathcal{O}(a),$$
implying $\Z_t(\{0\})=0$ for all $t\geq 0$.
Now putting $x=0$ in \eqref{theo2}, we have\footnote{corrected with - below, last line}
\begin{eqnarray*}
\Z_t(C_0) &=& F(S(0,t))+\int_0^tG^c(S(s,t))dB_s
\\
&=&F(S(0,t))+\la\int_0^tG^c(S(s,t))ds-\int_0^tG^c(S(s,t))dQ_s
\\
&=& F(S(0,t))+\la\int_0^tG^c(S(s,t))ds-Q_tG^c(0)+Q_0G^c(S(0,t))-\int_0^tg(S(s,t))Q_s d(S(s,t)),
\end{eqnarray*}
where we used integration-by-parts formula in the last line. Since $G^c(0)=1$, using \eqref{theo1} we have
\begin{eqnarray*}
X_t =Z_t+Q_t &=&F(S(0,t))+Q_0G^c(S(0,t))+\la\int_0^tG^c(S(s,t))ds
\\
&&\ -\int_0^t(X_s-1)^+  g(S(s,t))d(S(s,t)).
\end{eqnarray*}
By our assumptions on $F(\cdot), G(\cdot)$ and $k(\cdot)$, we see that $S(\cdot,t)$ is Lipschitz continuous on $[0,t]$ and hence $X_t$ is continuous in $t$.
Therefore
\begin{eqnarray*}
X_t =Z_t+Q_t &=&F(S(0,t))+Q_0G^c(S(0,t))+\la\int_0^tG^c(S(s,t))ds
\\
&&\ +\int_0^t(X_s-1)^+  k(X_s)g(S(s,t))ds.
\end{eqnarray*}
Hence by \eqref{theo1},
$Z_t, Q_t$ are continuous in $t$.
Now defining $H_1(x)=F(x)+Q_0G^c(x), H_2(x)=\la G^c(x),  H_4(x)=g(x), H_5(x)=k(x)$ for $x\geq 0$ and extending these maps on $(-\infty,0]$ by their
respective values at $0$, we have
\be
X_t =H_1(S(0,t))+\int_0^tH_2(S(s,t))ds
+\int_0^t(X_s-1)^+ H_4(X_s) H_5(S(s,t)) ds.\label{theo3}
\ee
By Lemma \ref{appn-lem2} in Appendix A, $X_t$ is uniquely defined on $[0, T]$ for all $T>0$. Since $Q_t=[X_t-1]^+$, $Q_t$ and hence $B_t$ is unique. Therefore
\eqref{theo2} implies that $\Z_t$ is uniquely defined on $\bar C_x$. Since $\{\bar C_x, x\in\R_+\}$ defines uniquely a Borel measure on $(\R_+,\BB)$, $\Z_t$ is uniquely defined by
\eqref{theo2}. Hence $(\Z, Q)$ is unique in $D([0,\infty),\M\times\R_+)$.\hfill $\Box$

\begin{remark}
One can relax the conditions on $G(\cdot)$ depending on the properties of service rate $k(\cdot)$. For example, if
$k(\cdot)$ is constant then it is enough to impose continuity on $G(\cdot)$ instead of Condition \ref{cond-1}(b)
(see \cite{zhang-dai-zwart, zhang}).
\end{remark}

\section{Tightness of the pre-limit processes}

In this section, we study the compactness properties of the pre-limit processes. From \eqref{4}, we get the following
equation
\be\label{3.1}
\bar\Z^n_t(C)=\bar\Z^n_s(C+S^n(s, t))+\frac{1}{n}\sum_{i=B^n_s+1}^{B^n_t}\del_{v^n_i}(C+S^n(\tau^n_i, t)),
\ee
for all $0\leq s\leq t$ and $C\in\BB$. Define $\bar E^n(s,t)=\bar E^n_t-\bar E^n_s$. From \eqref{7}, it is easy to see that
given $T, \eps>0$,
\be\label{3.2}
\lim_{n\to\infty}\p_n(\sup_{0\leq s\leq t\leq T}|\bar E^n(s,t)-\la(t-s)|\leq\eps)=1.
\ee
Defining $\bar B^n_t=\frac{B^n_t}{n},$ we have
\be\label{3.13}
\bar B^n_t=\bar E^n_t-\bar Q^n_t.
\ee
We recall the following characterization of compact subsets of $\M$ in Prohorov topology from \cite{delay-verejones} (Theorem A2.4.I).
\begin{definition}\label{def3.1}
A set ${\bf K}\subset\M$ is relatively compact if and only if $\sup_{\mu\in {\bf K}}\mu(\R_+)<\infty$ and there exists a sequence of nested compact
sets $C_j\subset\R_+$ such that $\cup C_j=\R_+$ and
$$ \lim_{j\to\infty}\sup_{\mu\in{\bf K}}\mu(C_j^c)=0.$$
\end{definition}
The proof of the following lemma is same as  \cite[Lemma 5.1]{zhang-dai-zwart}.
\begin{lemma}\label{lem1}
Fix $T>0$. There exists a sequence $\{\eps_E(n)\}$ such that $\eps_E(n)\to 0$ as $n\to\infty$ and
$$ \p_n(\sup_{0\leq s\leq t\leq T}|\bar E^n(s,t)-\la(t-s)|\leq\eps_E(n))\geq 1-\eps_E(n),$$
for all $n\in\N$.
\end{lemma}
We define $\Om^n_E=\{\sup_{0\leq s\leq t\leq T}|\bar E^n(s,t)-\la(t-s)|\leq\eps_E(n)\}$. The following lemma proves compact containment
properties of $(\bar\Z^n, \bar Q^n)$.

\begin{lemma}\label{lem2}
Assume Condition \ref{cond-1} to hold. Fix $T>0$. Then for any positive $\eta$ there exists a compact set ${\bf K}\subset \M$ and $K>0$ such that
$$\liminf_{n\to\infty}\p_n(\bar\Z^n_t\in{\bf K} \ \mbox{and}\ \bar Q^n_t\leq K\ \mbox{for all}\ t\in [0, T])\geq 1-\eta.$$
\end{lemma}

\noi{\bf Proof:}\ By Condition \ref{cond-1}(c), there exists a positive integer $M_0$ such that
$$ \sup_n\p_n(\bar Q^n_0>M_0)<\frac{\eta}{8}.$$
Since $\bar Q^n_t\leq \bar Q^n_0+\bar E^n_t$, choosing $K=M_0+2\la T$ we get from \eqref{3.2} that
\be\label{3.3}
\limsup_{n\to\infty}\p_n(\bar Q^n_t>K \ \mbox{for all}\ t\in[0, T])<\frac{\eta}{4}.
\ee
Define $\vartheta^n_t= \frac{1}{n}\sum_{i=-Q^n_0+1}^{E^n_t}\del_{v^n_i}\in \M$.
From \eqref{3.1}, we note that
\be\label{3.8}
\bar\Z^n_t(C_x)\leq \bar\Z^n_0(C_x)+\vartheta^n_t(C_x),
\ee
for all $x\in\R_+$ and $t\geq 0$. For $m\in\mathbb{Z}, \ell\geq 0$, define $\Li^n(m, \ell)=\frac{1}{n}\sum_{i=1+m}^{m+\lfloor n\ell\rfloor}\del_{v^n_i}$. Recall the definition of $\Om_A^n(M, L)$ (see \eqref{a2} in Appendix) and function $\bar f$ from Appendix B.
Define $\Om_1^n=\{\bar Q^n_0\leq M_0\}$. Therefore using Lemma \ref{lem1} we have for all large $n$,
 $$\p_n(\Om^n_E\cap\Om^n_1\cap\Om_A^n(M_0+1, K))\geq 1-\frac{\eta}{4}.$$
Choose $n$ large enough so that $\eps_A(n)\leq 1$. Then on $\Om^n_E\cap\Om^n_1\cap\Om_A^n(M_0+1, K)$ we have for all large $n$
\be\label{3.4}
\lan \bar f, \vartheta^n_t\ran\leq  K \lan \bar f, \nu^n\ran +1\leq \lan \bar f, \bar\nu\ran +1\leq M_1,
\ee
for some positive constant $M_1$. Hence using Markov's inequality and \eqref{3.4} we get
\be\label{3.5}
\vartheta^n_t(C_x)\leq \frac{1}{\bar f(x)}M_1,
\ee
on $\Om^n_E\cap\Om^n_1\cap\Om_A^n(M_0+1, K)$ for all $t\in[0, T]$ and all $n$ large. Again by Condition \ref{cond-1}(c),
there exists a compact ${\bf K}_0\subset\M$ such that for all large $n$
\be\label{3.6}
 \p_n(\bar\Z^n_0\in {\bf K}_0)\geq 1-\frac{\eta}{4}.
\ee
We denote the above event by $\Om^n_2$. By Definition \ref{def3.1}, there exists a bounded function $\varrho:\R_+\to\R_+$ such that $\lim_{x\to\infty}\varrho(x)=0$ and on $\Om^n_2$,
\be\label{3.7}
\bar\Z^n_0(\R_+)\leq \rho(0),\qquad \bar\Z^n_0(C_x)\leq\varrho(x)\ \qquad \mbox{for all}\ x\in\R_+.
\ee
We define
$${\bf K}=\{\mu\in\M : \mu(\R_+)\leq\rho(0)+K,\ \mu(C_x)\leq \varrho(x)+\frac{1}{\bar f(x)}M_1\ \forall \ x\in\R_+\}.$$
By the property of $\varrho$ and $\bar f$, ${\bf K}$ is a compact subset of $\M$. From \eqref{3.8}, \eqref{3.5} and
\eqref{3.7} we see that for all large $n$, on $\Om^n_E\cap\Om^n_1\cap\Om_A^n(M_0+1, K)\cap\Om^n_2$,
$$\bar\Z^n_t\in{\bf K},$$
for all $t\in[0, T]$. Also for all large $n$, $\p_n(\Om^n_E\cap\Om^n_1\cap\Om_A^n(M_0+1, K)\cap\Om^n_2)\geq 1-\frac{\eta}{2}$. Hence the proof follows from \eqref{3.3}.\hfill $\Box$

\begin{lemma}\label{lem3}
Assume Conditions \ref{cond-1}, \ref{cond-2} to hold. Fix $T>0$. Then for each $\eps, \eta>0,$ there exists a $\kappa>0$
such that
\be\label{3.9}
\liminf_{n\to\infty}\p_n(\sup_{[0,T]}\sup_{x\in\R_+}\bar\Z^n_t([x, x+\kappa])\leq\eps)> 1-\eta.
\ee
\end{lemma}

\noi{\bf Proof:}\ Using Condition \ref{cond-1}(c), one can prove that for any $\eps, \eta>0$ there
exists a positive $\kappa$ such that
\be\label{3.10}
\liminf_{n\to\infty}\p_n(\sup_{x\in\R_+}\bar\Z^n_0([x, x+\kappa])\leq\eps/2)\geq 1-\frac{\eta}{4}.
\ee
In fact, the proof is same as the proof of (79) in \cite{zhang-dai-zwart}. We denote the above event by $\Om^n_3$
and the event in Lemma \ref{lem2} by $\Om^n_4$. Define $\Om^n_5=\Om^n_3\cap\Om^n_4\cap\Om^n_E\cap\Om^n_A(M, L)$ for
$L=2\la T+K$ and $M=\lfloor L \rfloor +1$. From Lemma \ref{lem2}, it is easy to see that there exists $K>0$ such that
$$\liminf_{n\to\infty}\p_n(\Om^n_5)\geq 1-\eta.$$
From \eqref{3.1}, we see that for any $\kappa>0$ and $t\in[0, T]$
\be\label{3.11}
\bar\Z^n_t([x, x+\kappa])=\bar\Z^n_0([x,x+\kappa]+S^n(0,t))+\frac{1}{n}\sum_{i=-Q^n_0+1}^{B^n_t}\del_{v^n_i}([x, x+\kappa]+S^n(\tau^n_i, t)).
\ee
On $\Om^n_3$, we have $\sup_{x\in\R_+}\bar\Z^n_0([x, x+\kappa])\leq\eps/2$. Hence choosing $x=x(\omega)=x+S^n(0,t)$
we have on $\Om^n_3$, $\sup_{x\in\R_+}\bar\Z^n_0([x,x+\kappa]+S^n(0,t))\leq\eps/2$. Thus we need to estimate the second term on
$\Om^n_5$. So we denote the second term by $\Xi_t$.

For any $\del>0$, we consider a partition $0=t_0<t_1<\ldots<t_r=t$ of $[0,t]$ with $|t_{i+1}-t_i|<\del$ for
$k=0,1,2,\ldots,r-1$. Since $-Q^n_0=B^n_0$, we have
$$ \Xi_t=\sum_{k=0}^{r-1}\frac{1}{n}\sum_{i=B^n_{t_k}+1}^{B^n_{t_{k+1}}}\del_{v^n_i}([x, x+\kappa]+S^n(\tau^n_i, t)).$$
From \eqref{a2}, on $\Om^n_A(M,L)$, we have for all $k=0,1,2,\ldots,r-1$,
\begin{displaymath}
-Q^n_0 \leq  B^n_{t_i}\leq E^n_t,\ \mbox{and so},\ \ 0\leq B^n_{t_{i+1}}-B^n_{t_i}\leq E^n_t+Q^n_0,
\end{displaymath}
and
\begin{displaymath}
\max_{-nM<m<nM}\sup_{\ell\in [0, L]}\sup_{f\in\mathcal{V}}|\lan f, \Li^n(m, \ell)\ran-\ell\lan f, \nu^n\ran|\leq\eps_A(n).
\end{displaymath}
Hence for all $m\in (-nM, nM), \ell\in [0,L]$ and for all $a, b\in\R_+, a\leq b,$ we have
\be\label{3.12}
\lan \chi_{[a,b]}, \Li(m, \ell)\ran\leq \ell\lan \chi_{[a,b]}, \nu^n\ran + 2\eps_A(n),
\ee
on $\Om^n_A(M,L)$. Now for $t_{k}\leq \tau^n_i\leq t_{k+1}$, $[x, x+\kappa]+S^n(\tau^n_i, t)\subset [x +S^n(t_{k+1}, t), x+\kappa+S^n(t_{k}, t)]$.
Now fixing $a=x +S^n(t_{k+1}, t), b=x+\kappa+S^n(t_{k}, t)$ and observing that, on $\Om^n_4\cap\Om^n_E$, $B_{t_i}\in(-nM, nM)$ (for above choice of $M$)
and $\bar B^n_{t_{i+1}}-\bar B^n_{t_i}\in [0, L]$ for all $n$ large ,
we have
$$\frac{1}{n}\sum_{i=B^n_{t_i}+1}^{B^n_{t_{i+1}}}\del_{v^n_i}([x, x+\kappa]+S^n(\tau^n_i, t))\leq (\bar B^n_{t_{i+1}}-\bar B^n_{t_i})
\nu^n([x +S^n(t_{k+1}, t), x+\kappa+S^n(t_{k}, t)])+2\eps_A(n),$$
for $k=0,1,2,\ldots,r-1$. Since $\nu^n\to\nu$ in Prohorov metric, for any $\eps_1>0$, there exists $n_0\in\N$ such that for all
$n\geq n_0$, and closed set $C\subset\R_+$ (see Notations)
$$\nu^n(C)\leq \nu(C^{\eps_1})+\eps_1.$$
Hence combining above two we have
\begin{eqnarray*}
\frac{1}{n}\sum_{i=B^n_{t_i}+1}^{B^n_{t_{i+1}}}\del_{v^n_i}([x, x+\kappa]+S^n(\tau^n_i, t)) &\leq &(\bar B^n_{t_{i+1}}-\bar B^n_{t_i})
\nu([x +S^n(t_{k+1}, t)-\eps_1, x+\kappa+S^n(t_{k}, t)+\eps_1])
\\
&&+(\bar B^n_{t_{i+1}}-\bar B^n_{t_i})\eps_1+2\eps_A(n),
\end{eqnarray*}
on $\Om^n_4\cap\Om^n_E$ for all $n$ large. At this point, we observe that $\sup_{s\in [0, T]}\bar X^n_s\leq \sup_{s\in[0,T]}\bar Q^n_s+1\leq L+1$
on $\Om^n_4$ for all large $n$. Hence by Condition \ref{cond-2}, on $\Om^n_4$, $\sup_{[0,T]}k^n(\bar X^n_s)\leq \sup_{[0,L+1]}k^n(x)\leq
\sup_{[0,L+1]}k(x)+1<M_2$ for some positive constant $M_2$ and for all $n$ large. Hence $|S^n(t_k, t)-S^n(t_{k+1}, t)|\leq\del M_2$ on
$\Om^n_4$ for all $n$ large. By Condition \ref{cond-1}(b), we can choose $\del$ and $\kappa$ small enough so that on $\Om^n_4$
$$ \nu([x +S^n(t_{k+1}, t)-\eps_1, x+\kappa+S^n(t_{k}, t)+\eps_1])<\frac{\eps}{8L},$$
for all $\eps_1$ small enough and all $n$ large. Hence summing up the above expression we have for all $t\in [0, T]$
$$ \Xi_t\leq L. \frac{\eps}{8L}+ L\eps_1+ 2r\eps_A(n),$$
on $\Om^n_5$ for all $n$ large. Since $\eps_1, \eps_A(n)$ do not depend on $r$ and $x$, we can choose them small to make the right hand side
smaller than $\frac{\eps}{2}$ for all $n$ large and $x\in\R_+, t\in[0,T]$. The proof is done from \eqref{3.11} and definition of $\Om^n_5$.\hfill $\Box$

\spa

For any path $\phi\in D([0, \infty), E)$ where $(E, \pi)$ is polish space, the $\del$-oscillation of $\phi$ on $[0, T],\ T>0$, is defined as follows
$$ w(\phi, \del, T)=\sup_{s, t\in[0,T], |s-t|\leq\del}\pi(\phi(s),\phi(t)).$$
The following lemma gives the oscillation bounds on the stochastic process
\begin{lemma}\label{lem4}
Assume Conditions \ref{cond-1}, \ref{cond-2} to hold. Fix $T>0$. Then for each $\eps, \eta>0$, there exists $\del>0$ such that
\begin{eqnarray*}
\liminf_{n\to\infty}\p_n(w(\bar\Z^n,\del, T)\leq\eps) &\geq  & 1-\eta,
\\
\liminf_{n\to\infty}\p_n(w(\bar Q^n,\del, T)\leq\eps) &\geq &1-\eta.
\end{eqnarray*}
\end{lemma}

\noi{\bf Proof:}\  Let $t, s\in[0, T], s\leq t$ and $|t-s|\leq\del$. Let $D^n_t$ be the
number of customers finished their job by time $t$. Then it is easy to see that
$$ D^n_t-D^n_s\leq \Z^n_s[0,S^n(s,t)]+\sum_{i=B^n_s+1}^{B^n_t}\del_{v^n_i}([0,S^n(s,t)]).$$
Recall the event $\Om^n_5$ from Lemma \ref{lem3}. By definition, on $\Om^n_5$, we have $B^n_s\in (-nM, nM)$ and $\bar B^n_{t}-\bar B^n_s\in
[0,L]$ for all $n$ large. Since $\sup_{[0, T]}\bar X^n_s\leq\sup_{[0,T]}\bar Q^n_s+1\leq L+1$, we can choose $\del$ small enough so that
$|S^n(s,t)|\leq \kappa_1$ on $\Om^n_5$ for all $n$ large where $\kappa_1=\kappa_1(\del)\to 0$ as $\del\to 0$. Therefore by the
definition of $\Om^n_A(M,L)$ and \eqref{3.12} we have
$$\frac{1}{n}\sum_{i=B^n_s+1}^{B^n_t}\del_{v^n_i}([0,S^n(s,t)])\leq L\lan \chi_{[0, \kappa_1]}, \nu^n\ran + 2\eps_A(n),$$
on $\Om^n_5$ for all $n$ large. Since $\nu^n\to\nu$, by the same reasoning as in Lemma \ref{lem3}, we can choose $\del$ small so that
$$\frac{1}{n}\sum_{i=B^n_s+1}^{B^n_t}\del_{v^n_i}([0,S^n(s,t)])\leq\eps,$$
on $\Om^n_5$ for all $n$ large. If we denote the event in \eqref{3.9} by $\Om^n_6$, then on $\Om^n_6\cap\Om^n_5$
$$ \bar\Z^n_s[0,S^n(s,t)]\leq \ \bar\Z^n_s[0,\kappa_1]\leq\eps$$
for all $n$ large and $\del$ chosen small enough. Hence with this choice of $\del$, $\p_n(\Om^n_6\cap\Om^n_5)\geq 1-2\eta$ and
$$ \frac{1}{n}(D^n_t-D^n_s)\leq 2\eps$$
for all $n$ large. Since $X^n_t=X^n_0+E^n_t-D^n_t$, we have
$|\bar X^n_s-\bar X^n_t|\leq |\bar E^n_s-\bar E^n_t|+ \frac{1}{n}(D^n_t-D^n_s)\leq 3\eps$ on $\Om^n_6\cap\Om^n_5$ for all $n$ large
provided $\del$ chosen small enough. Therefore with this choice of $\del$, we have (using \eqref{6})
$$|\bar Q^n_s-\bar Q^n_t|\leq3\eps$$
on $\Om^n_6\cap\Om^n_5$ for all $n$ large. Hence the second claim follows by replacing $\eps, \eta$ with $\eps/3, \eta/2$ respectively.

Now we prove the first claim. We note that for any $t,s\in[0,T]$, $|\bar B^n_s-\bar B^n_t|\leq |\bar E^n_s-\bar E^n_t|+|\bar Q^n_s-\bar Q^n_t|$
(from \eqref{3.13}). Hence on $\Om^n_6\cap\Om^n_5$, $|\bar B^n_s-\bar B^n_t|\leq 4\eps$ for all $n$ large provided $|t-s|\leq\del\wedge\frac{\eps}{2\la}$. Let
$C\subset\R^+$ be closed and $C^\eps$ be its $\eps$-enlargement. Now choose $\del$ small enough so that $\kappa_1<\eps$ and so $C+S^n(s,t)\subset C^{\kappa_1}\subset C^\eps$ on $\Om^n_5$ for all $n$ large. Hence from \eqref{3.1} we have
$$ \bar\Z^n_t(C)-\bar\Z^n_s(C^\eps)\leq \bar\Z^n_s(C+S^n(s,t))-\bar\Z^n_s(C^\eps)+ |\bar B^n_s-\bar B^n_t|\leq 4\eps,$$
on $\Om^n_6\cap\Om^n_5$ for all $n$ large.
Again for any $c\in C$ and $S^n(s,t)\leq \kappa_1$, we have $\text{dist}(c-S^n(s,t), C)<\eps$ implying $c\in C^\eps+S^n(s,t)$
and so $C\subset C^\eps+S^n(s,t)$. Hence from \eqref{3.1}
$$  \bar\Z^n_s(C)-\bar\Z^n_t(C^\eps)\leq \bar\Z^n_s(C)-\bar\Z^n_s(C^\eps+S^n(s,t))\leq 0,$$
on $\Om^n_6\cap\Om^n_5$ for all $n$ large. Hence for all closed set $C\in\BB$ we have
$$ \bar\Z^n_s(C)\leq \bar\Z^n_t(C^{4\eps})+4\eps\ \mbox{and}\ \bar\Z^n_t(C)\leq \bar\Z^n_s(C^{4\eps})+4\eps,$$
on $\Om^n_6\cap\Om^n_5$ for all $n$ large. Hence $\rho(\bar\Z^n_s,\bar\Z^n_t)\leq 4\eps$ on $\Om^n_6\cap\Om^n_5$ for all $n$ large. Thus
first claim follows by replacing $\eps, \eta$ with $\eps/4, \eta/2$ respectively. \hfill $\Box$

\spa

\noi Now we introduce a weaker oscillation function $w^\prime$ on $D([0,\infty),\M\times\R_+)$. Define metric $d'(\cdot, \cdot)=\max\{\rho(\cdot,\cdot),|\cdot|\}$
on $\M\times\R_+$ which induces a separable complete metric on it. For any $\psi\in D([0,\infty),\M\times\R_+)$ and $T, \del>0$ define
$$ w^\prime(\psi, \del, T)=\inf_{t_i}\max_{i}\sup_{s,t\in[t_{i-1}, t_i)}d'(\psi(s),\psi(t)),$$
where $\{t_i\}$ ranges over all partition of the form $0=t_0<t_1<\ldots<t_j=T$ with $\min_{1\leq i\leq j}(t_i-t_{i-1})>\del$ and $j\geq 1$.
It is easy to see that for $\del>0$ we can have a partition ${t_i}$ of $[0, T]$ such that $\min_{1\leq i\leq j}(t_i-t_{i-1})>\del$
and $\max_{1\leq i\leq j}(t_i-t_{i-1})\leq 2\del$ and hence for any $\psi\in D([0,\infty),\M\times\R_+)$ we have
\be\label{3.17}
w^\prime(\psi, \del, T)\leq w(\psi, 2\del, T).
\ee
Hence from Lemma \ref{lem4} and \eqref{3.17}, we get that for any $T, \eps, \eta>0,$ there exists $\del>0$ such that
\be\label{3.14}
\liminf_{n\to\infty}\p_n(w'((\bar\Z^n, \bar Q^n),\del, T)\leq\eps) \geq   1-\eta.
\ee
For $\psi\in D([0,\infty),\M\times\R_+)$, define
$$J(\psi)=\int_0^\infty e^{-s}[J(\psi, s)\wedge 1]ds, \ \mbox{where}\ \ J(\psi, t)=\sup_{0\leq s\leq t}d'(\psi(s),\psi(s-)).$$
Again it is easy to see that $J(\psi, t)\leq w(\psi, \del, t)$ for any $\del>0$ and hence $J(\psi)\leq w(\psi,\del, T)+e^{-T}$
for all $T, \del>0$. Thus applying Lemma \ref{lem4}, we get that for any $\eps, \eta>0,$
\be\label{3.15}
\liminf_{n\to\infty}\p_n(J((\bar\Z^n,\bar Q^n))\leq\eps)\geq 1-\eta.
\ee
Now we note that the process $(\bar\Z^n,\bar Q^n)$ satisfies the conditions (a) compact containment property (Lemma \ref{lem2})
and (b) oscillation bound (\eqref{3.14}) of Corollary 3.7.4 in \cite{ethier-kurtz} . Hence $(\bar\Z^n,\bar Q^n)$
is tight in $D([0,\infty),\M\times\R_+)$ and $(\bar\Z^n,\bar Q^n)\Rightarrow (\bar\Z,\bar Q)$ along some subsequence for some
random variable $(\bar\Z,\bar Q)$ taking values in $D([0,\infty),\M\times\R_+)$. Also \eqref{3.15} satisfies the condition (a) in Theorem 3.10.2 in \cite{ethier-kurtz} which implies that $(\Z, Q)$ has continuous paths almost surely.

\section{Characterization of the limits}

In this section, we characterize some properties of the limits which lead to uniqueness. To have simple notations, we consider
full sequence to converge instead of subsequence.
To this end, we intend to define all the variable
on a common probability space using Skorohod representation theorem.
From Condition \ref{cond-1}(a), it is clear that $\bar E^n$ is tight in $D([0,\infty),\R_+)$. Hence
$(\bar\Z^n,\bar Q^n,\bar E^n)$ is tight in $D([0,\infty),\M\times\R_+)\times D([0,\infty),\R_+)$. Therefore by Skorohod representation theorem we can say that
${\color{red}(\bar\Z^{1n},\bar Q^{1n},\bar E^{1n}, U^{1n}, V^{1n})} \to (\Z^1,Q^1,\la\cdot)$ in $D({\color{red}[0,T]},\M\times\R_+)\times D([0,\infty),\R_+)$ almost surely on some probability space $(\tilde\Om, \tilde{\mathcal{F}},\tilde\p)$ where
\begin{center}
law of $(\bar\Z^{1n},\bar Q^{1n},\bar E^{1n})$ = law of $(\bar\Z^n,\bar Q^n, \bar E^n)$ for all $n$,
\end{center}
and law of $(\Z^{1}, Q^{1})$ = law of $(\Z, Q)$.

Also $(Z^1, Q^1)$ has continuous paths almost surely for $Z^1_t=\lan 1, \Z^1_t\ran$. Hence for any $T>0$, the followings hold, almost surely:
\begin{eqnarray}
\lim_{n\to\infty}\sup_{s\in [0, T]}\rho(\Z^{1n}_s,\Z^1_s) &=& 0\label{4.8}
\\
\lim_{n\to\infty}\sup_{s\in[0, T]}|\bar Q^{1n}_s-Q^1_s| &= & 0,\nonumber
\\
\lim_{n\to\infty}\sup_{s\in[0, T]}|\bar Z^{1n}_s-Z^1_s| &=& 0,\nonumber
\\
\lim_{n\to\infty}\sup_{s\in[0, T]}|\bar B^{1n}_s-B^1_s| &=& 0,\label{4.1}
\end{eqnarray}
where $Z^1_t=\lan 1, \Z^1_t\ran$, $B^{1n}_t=E^{1n}_1-Q^{1n}_t$ and $B^1_t=\la t-Q^1_t$. Hence from \eqref{5} and \eqref{6}, we get $\bar X^{1n}\to X^1$ uniformly on $[0,T], T>0,$ and
$$ X^1_t=Q^1_t+Z^1_t\ \mbox{and}\ Q^1_t=[X^1_t-1]^+.$$
\eqref{4.1} implies that $B^1_t$ is nondecreasing in $t$ and so it is a function of bounded variation on $[0,T]$ for all $T>0$. Also from above, it is easy to see that
\begin{center}
law of $(\Z^1, Q^1, B^1, X^1)$= law of $(\Z, Q, B, X)$.
\end{center}

\begin{lemma}\label{emni}
Let $(\mathcal{S},\pi)$ be a metric space and $K\subset\Sp$ be compact. Let $f:\Sp\to\R$ be a function satisfying
the following: for any sequence $s_n\to s$ and $s\in K$, $f(s_n)\to f(s)$ as $n\to\infty$. Then for any $\eps>0$,
there exists a $\del>0$ such that $|f(s_1)-f(s_2)|<\eps$ whenever $s_1\in\Sp, s_2\in K$ and
$\pi(s_1,s_2)\leq\del$.
\end{lemma}

\noi{\bf Proof:}\ If not, then there exists two sequences $\{s_n\},\{\tilde s_n\}$ such that $\{\tilde s_n\}\subset K$ and
$$ \pi(s_n,\tilde s_n)\leq\frac{1}{n}, \ |f(s_n)-f(\tilde s_n)|\geq\eps \ \ \forall\ n\geq 1.$$ Now $K$ being compact,
there exists $s\in K$ such that along some subsequence $\{n_k\}$, $\tilde s_{n_k}\to s\in K$. Hence $s_{n_k}\to s\in K$
as $n_k\to\infty$. This is contradicting to the fact that $|f(s_{n_k})-f(\tilde s_{n_k})|\geq\eps$ for all $n_k$.
Hence the proof.\hfill $\Box$

\spa

Following lemma is a consequence of Lemma \ref{lem3}.

\begin{lemma}\label{lem4.6}
Fix $T>0$ and $x_0\in\R_+$. Consider a decreasing sequence $\{f^n\}$ in $ C_b(\R)$ so that $f^n\geq0, f^n(x)=1$ on $[x_0-\frac{1}{n}, x_0+\frac{1}{n}]$ and $f^n$ vanishes outside of $[x_0-\frac{2}{n}, x_0+\frac{2}{n}]$.
Then $\tilde\p(\lim_{n\to\infty}\sup_{t\in[0,T]}\langle f^n, \Z^1\rangle>0)=0$. In particular, for $t\geq 0$, $\Z^1_t$ has
no atom at $x_0$ almost surely.
\end{lemma}

\noi{\bf Proof:}
 Let $\tilde\p(\lim_{n\to\infty}\sup_{t\in[0,T]}\langle f^n, \Z^1\rangle>\kappa_2)\geq \eta$ for some positive constant
$\kappa_2, \eta$. Then $\tilde\p(\sup_{t\in[0,T]}\lan f^m, \Z^1_t\ran>\kappa_2)\geq \eta$ for all $m$.
Note that $\{\Z^1_t : t\in[0,T]\}$ is compact in $\M$.
Now from \eqref{4.8} and Lemma \ref{emni} we get, $\sup_{t\in[0,T]}\lan f^m,\bar\Z^{1,n}\ran\to \sup_{t\in[0,T]}\lan f^m,\Z^1\ran$ as $n\to\infty$ almost surely.
Therefor using Fatou's lemma, for all $m$,
\begin{eqnarray*}
 \liminf\p_n(\sup_{t\in[0,T]}\lan f^m,\bar \Z^n_t\ran>\kappa_2) =\liminf \tilde\p(\sup_{t\in[0,T]}\lan f^m,\bar \Z^{1n}_t\ran>\kappa_2)\geq\tilde\p(\sup_{t\in[0,T]}\lan f^m, \Z^1_t\ran>\kappa_2)\geq\eta.
\end{eqnarray*}
Now we choose $\kappa>0$ from Lemma \ref{lem3} for $\eps, \eta$ replaced by $\kappa_2/2, \eta/2$.
Thus if we choose $m$ large enough, we get
$$\eta\leq \liminf_{n\to\infty}\p_n(\sup_{t\in[0,T]}\lan f^m,\bar \Z^n_t\ran>\kappa_2)\leq
\liminf_{n\to\infty}\p_n(\sup_{t\in[0,T]}\bar \Z^n_t[(x_0-\frac{\kappa}{4})\vee 0, x_0+\frac{\kappa}{4}]>\kappa_2)\leq \frac{\eta}{2}, $$
which is a contradiction. This completes the proof.\hfill $\Box$.

An immediate consequence of the above lemma is $\Z_t([0,\infty)=\Z_t((0,\infty))=Z_t$
for all $t\geq 0,$ almost surely.

\begin{lemma}
For any $t\geq 0$, $\Z_t$ satisfies the fluid model equation $(k, \la, \nu)$ given by
$$\Z_t(\bar C_x)=\Z_0(\bar C_x+S(0,t))+\int_0^tG^c(x+S(s,t))dB_s,$$
almost surely where $S(s,t)=\int_s^tk(X_u)du$.
\end{lemma}

\noi{\bf Proof:}\ It is enough to prove the above result for the process $\Z^1$.
For $t\geq 0$ and $C\in\BB$, we have
\be\label{4.2}
\Z^n_t(C)=\Z_0(C+S^n(0,t))+I^n_t(C),
\ee
where
$$I^n_t(C)=\sum_{j=0}^{J-1}\frac{1}{n}\sum_{i=B^n_{t_j}+1}^{B^n_{t_{j+1}}}\del_{v^n_i}(C+S^n(\tau^n_i,t)),$$
for any partition $\{t_j\}_{j=0}^{j=J(t)}$ of $[0,t]$.
But $\Z^{1n}_t$ might not possess the same expression as $\Z^n_t$ as the stochastic variables ${v^i_n}$ might not make sense on new probability space
$(\tilde\Om,\tilde{\mathcal{F}},\tilde\p)$. Define $S^{1n}(s,t)=\int_s^tk^n(\bar X^{1n}_u)du$
and $S^{1}(s,t)=\int_s^tk( X^{1}_u)du$. Fix $x\in\R_+$ and let $A_T$ be a countable dense set in $[0, T]$. Then for any $\eps>0$,
\begin{eqnarray}
&&\tilde\p(\sup_{t\in A_T}|\Z^{1n}_t(\bar C_x)-\Z^{1n}_0(\bar C_x+S^{1n}(0,t))-\int_0^tG^c(x+S^{1n}(s,t))d\bar B^{1n}_s|>\eps)\nonumber
\\
&&=\p_n(\sup_{t\in A_T}|\Z^{n}_t(\bar C_x)-\Z^{n}_0(\bar C_x+S^n(0,t))-\int_0^tG^c(x+S^n(s,t))d\bar B^{n}_s|>\eps)\nonumber
\\
&&=
\p_n(\sup_{t\in[0,T]}|I^n_t(\bar C_x)-\int_0^tG^c(x+S^n(s,t))d\bar B^n_s|>\eps).\label{4.3}
\end{eqnarray}
Applying Lemma \ref{lem2}, for any positive $\eta$, we have constant $K$ such that $\p_n(\Om^n_7)\geq 1-\eta$ where $\Om^n_7=\{\sup_{t\in [0,T]}\bar X^n_t\leq K\}$
for all $n$ large.


We choose $\del>0$ and partitions $\{t_i\}_{i=0}^{i=J(t)}$ such that
$\max_{1\leq j\leq J(t)}|t_{j}-t_{j-1}|<\del$ and $\sup_{t\in[0,T]}J(t)=J(\del)<\infty$.
Recall that $V^n\to 0$ as $n\to\infty$ in probability. From Lemma \ref{lem1} and Lemma \ref{lem2}, we choose $M$ such that
$\p_n(\Om^n_8)\geq 1-\eta$ for all $n$ large where $\Om^n_8=\{\sup_{[0,T]}\bar B^n_s< \frac{M}{2}\}$.
Hence on $\Om^n_8$, $(\bar B^n_{t_{j+1}}-\bar B^n_{t_j})< M$ for all $n$.
Again for $i\in \{B^n_{t_j}+1, B^n_{t_j}+2,\ldots, B^n_{t_{j+1}}\}$, we have $t_j< \tau^n_i\leq t_{j+1}$. Recall the event $\Om^n_A(M, M)$ from \eqref{a2}.
For any $\eps_1>0$, we have for $0\leq j\leq J(t)-1,$
\begin{eqnarray*}
&&\frac{1}{n}\sum_{i=B^n_{t_j}+1}^{B^n_{t_{j+1}}}\del_{v^n_i}(\bar C_{x+S^n(\tau^n_i,t)})-\frac{1}{n}\sum_{t_j<\tau^n_i\leq t_{j+1}}G^c(x+S^n(\tau^n_i,t))
\\
&\leq &\frac{1}{n}\sum_{i=B^n_{t_j}+1}^{B^n_{t_{j+1}}}\del_{v^n_i}(\bar C_{x+S^n(t_{j+1},t)})-(\bar B^n_{t_{j+1}}-\bar B^n_{t_j})G^c(x+S^n(t_{j},t))
\\
&\leq & (\bar B^n_{t_{j+1}}-\bar B^n_{t_j})\Big(\nu^n(\bar C_{x+S^n(t_{j+1},t)})-G^c(x+S^n(t_{j},t))\Big)+ \eps_1
\end{eqnarray*}
on $\Om^n_A(M,M)\cap\Om^n_8$, for all $n$ large. Using the fact that $\rho(\nu^n,\nu)\to 0$ as $n\to\infty$ , we get on $\Om^n_A(M, M)\cap\Om^n_8\cap\Om^n_7$,
\begin{eqnarray*}
&&\sum_{j=0}^{J-1}\frac{1}{n}\sum_{i=B^n_{t_i}+1}^{B^n_{t_{i+1}}}\del_{v^n_i}(\bar C_{x+S^n(\tau^n_i,t)})-\int_0^tG^c(x+S^n(s,t))d\bar B^{n}_s
\\
&\leq &\sum_{j=0}^{J-1}(\bar B^n_{t_{j+1}}-B^n_{t_j})\Big(\nu(\bar C_x+(S(t_{j+1},t)-\eps_1)\vee 0)-G^c(x+S^n(t_{j},t))\Big)+ M\eps_1+ J(t)\eps_1,
\\
&\leq &M|g|_\infty(K\del+\eps_1)+ (M+J(\del))\eps_1,
\end{eqnarray*}
for all $n$ large where $|g|_\infty$ denote the supremum norm of $g$. First choosing $\del>0$ small enough and then choosing $\eps_1$ we can have the r.h.s.
less than $\eps/2$ on $\Om^n_A(M, M)\cap\Om^n_8\cap\Om^n_7$ for all $n$ large and for all $t\in[0,T]$. A similar calculation gives that $I^n_t(\bar C_x)-\int_0^tG^c(x+S^n(s,t))d\bar B^n_s\geq -\eps/2$
on $\Om^n_A(M, M)\cap\Om^n_8\cap\Om^n_7$ for all $n$ large  and for all $t\in[0,T]$. Since $\liminf_{n\to\infty}\p_n(\Om^n_A(M, M)\cap\Om^n_8\cap\Om^n_7)\geq 1-3\eta$, we have
from \eqref{4.3}
\begin{eqnarray*}
\limsup_{n\to\infty}\tilde\p(\sup_{t\in[0,T]}|\Z^{1n}_t(\bar C_x)-\Z^{1n}_0(\bar C_x+S^{1n}(0,t))-\int_0^tG^c(x+S^{1n}(s,t))d\bar B^{1n}_s|>\eps)\leq 3\eta.
\end{eqnarray*}
$\eta$ begin arbitrary, we have for any $\eps>0$,
\be\label{4.7}
\limsup_{n\to\infty}\tilde\p(\sup_{t\in[0,T]}|\Z^{1n}_t(\bar C_x)-\Z^{1n}_0(\bar C_x+S^{1n}(0,t))-\int_0^tG^c(x+S^{1n}(s,t))d\bar B^{1n}_s|>\eps)=0.
\ee

Since $\sup_{[0,T]}|\bar X^{1n}_s-X^1_s|\to 0$, by Condition \ref{cond-2}, $\sup_{[0,T]}|k^n(\bar X^{1n}_s)-k(X^1_s)|\to 0$
as $n\to\infty$, almost surely. Hence
\begin{eqnarray*}
\sup_{s\in[0,T]}|S^{1n}(s,t)-S^1(s,t)|\to 0\ \mbox{as}\ n\to\infty.
\end{eqnarray*}
This implies $\sup_{s\in[0,T]}|G^c(x+S^{1n}(s,t))-G^c(x+S^1(s,t))|\to 0$ as $n\to\infty$. Since $\rho(\Z^{1n}_0,\Z^1_0)\to 0$ almost surely, we have for $t\in [0,T]$
$$ \Z^{1n}_0(\bar C_x+S^{1n}(0,t))\leq \Z^{1}_0(\bar C_x+S^1(0,t)-\eps_2)+\eps_2, \Z^{1}_0(\bar C_x+S^1(0,t)+\eps_2)\leq \Z^{1n}_0(\bar C_x+S^{1n}(0,t))+\eps_2,$$
for any chosen $\eps_2>0$ and all $n$ large (might depend on sample point).
By Condition \ref{cond-1}(c), $\Z^1_0$ is deterministic with distribution function Lipschitz continuous and so $\sup_{t\in[0,T]}|\Z^{1n}_0(\bar C_x+S^{1n}(0,t))-\Z^{1}_0(\bar C_x+S^1(0,t))|\to 0$
almost surely. From \eqref{4.1}, it is easy to check that $\rho(d\bar B^{1n},dB^1)\to 0$ where $d\bar B^{1n}, dB^1$ are considered as Borel measures on $[0, T]$. Since $B^1$ is continuous almost surely,
applying Theorem A2.3.I in \cite{delay-verejones} and \eqref{4.1}, we get
$$ \sup_{t\in[0,T]}|\int_0^tG^c(x+S^1(s,t))d\bar B^{1n}_s-\int_0^tG^c(x+S^1(s,t))dB^1_s|\to 0\ \mbox{as}\ \ n\to\infty,\ \mbox{almost surely,}$$
and hence
$$\sup_{t\in[0,T]}| \int_0^tG^c(x+S^{1n}(s,t))d\bar B^{1n}_s-\int_0^tG^c(x+S^1(s,t))dB^1_s|\to 0 \ \mbox{as}\ \ n\to\infty,\ \mbox{almost surely.}$$
Now we show that $\sup_{t\in[0,T]}|\bar\Z^{1n}_t(\bar C_x)-\Z^1_t(\bar C_x)|\to 0$ as $n\to\infty$ almost surely.
Consider the map $f:\M\to\R$ defined by $f(\mu)=\mu(\bar C_x)$.
From Lemma \ref{lem4.6}, we see that $f$ satisfies the condition of Lemma \ref{emni} for the compact set
$\{\Z^1_t : t\in[0,T]\}$, almost surely. Hence using \eqref{4.8} and Lemma \ref{emni}, we have
$\sup_{t\in[0,T]}|\bar\Z^{1n}_t(\bar C_x)-\Z^1_t(\bar C_x)|\to 0$ as $n\to\infty$, almost surely.
Combining the above estimates with \eqref{4.7}, we get
$\Om(x,T)\in\tilde{\mathcal{F}}$ such that $\tilde\p(\Om(x,T))=1$ and
$$\Z^1_t(\bar C_x)=\Z^1_0(\bar C_x+S^1(0,t))+\int_0^tG^c(x+S^1(s,t))dB^1_s,$$
for all $t\in[0,T]$ on $\Om(x,T)$. Since $\{\bar C_x: x\in\R_+, x\ \mbox{rational}\}$ determines any Borel-measure uniquely on $\R_+$ we can take
$\Om_\infty=\cap_{T\in\N}\cap_{\{x: x\ rational\}}\Om(x,T)$ on which
$$\Z^1_t(\bar C_x)=\Z^1_0(\bar C_x+S^1(0,t))+\int_0^tG^c(x+S^1(s,t))dB^1_s,$$
for all $t\geq 0$ and $x\in\R_+$. This completes the proof as $\tilde\p(\Om_\infty)=1$.\hfill $\Box$

\section{Appendix}
\subsection{Appendix A}
In this section, we prove existence of unique solution to the fluid model type equations.
\begin{lemma}\label{appn-lem1}
 Fix $T>0$.
Let $k,H_i:\R\to\R, i=1,2,3,4,$ be bounded, Lipschitz continuous. Then the following integral equation
\begin{eqnarray}
x_t &=& H_1(\int_0^tk(x_s)ds)+\int_0^tH_2(\int_s^tk(x_u)du)ds+\int_0^tH_3(x_s)H_4(\int_s^tk(x_u)du)ds\label{appn-a1}
\\
x_0 &=& H_1(0).\nonumber
\end{eqnarray}
has a unique solution in $C([0,T], \R)$.
\end{lemma}

\noi{\bf Proof:}\ To simplify the notation, we define $S(s,t,\phi)=\int_s^tk(\phi(u))du$ for $\phi\in C([s,t],\R)$.
Assume that solution $x_t$ is uniquely defined on $[0,t_0]$ for $t_0\in[0,T)$. We consider the following integral
equation for $t\in [t_0, T]$
\begin{eqnarray}
x_t &=&H_1(S(0,t_0,x)+S(t_0, t, x))+ \int_0^{t_0}H_2(S(s,t_0,x)+S(t_0,t,x))ds\nonumber
\\
&&+\int_0^{t_0}H_3(x_s)H_4(S(s,t_0,x)+S(t_0,t,x))ds+ \int_{t_0}^tH_2(S(s,t,x))ds \nonumber
\\
&&+ \int_{t_0}^tH_3(x_s)H_4(S(s,t,x))ds.\label{apen-a2}
\end{eqnarray}
Now define a operator $F: C([t_0,T],\R)\to C([t_0,T],\R)$ as follows:
\begin{eqnarray*}
F(\phi)(t) &=&H_1(S(0,t_0,x)+S(t_0, t, \phi))+ \int_0^{t_0}H_2(S(s,t_0,x)+S(t_0,t,\phi))ds
\\
&&+\int_0^{t_0}H_3(x_s)H_4(S(s,t_0,x)+S(t_0,t,\phi))ds+ \int_{t_0}^tH_2(S(s,t,\phi))ds
\\
&&+ \int_{t_0}^tH_3(\phi(s))H_4(S(s,t,\phi))ds.
\end{eqnarray*}
To simplify the notation, we denote the $i$-th term on the r.h.s. of the above expression by $F_i(\phi)$ for $i=1,2,3,4,5$.
We denote the supremum (Lipschitz constant) of $H_i$ by $H_{i\infty}(L_i)$ for $i=1,2,3,4$. Let $L_k$ be the Lipschitz constant of $k(\cdot)$.
Then for $\phi^1, \phi^2\in C([t_0,T],\R)$ the followings hold: for $t\in [t_0, T]$
\begin{eqnarray*}
|F_1(\phi^1)(t)-F_1(\phi^2)(t)| &\leq & L_1L_k(t-t_0)|\phi^1-\phi^2|_{t_0t},
\\
|F_2(\phi^1)(t)-F_2(\phi^2)(t)| &\leq & L_2L_k t_0(t-t_0)|\phi^1-\phi^2|_{t_0t},
\\
|F_3(\phi^1)(t)-F_3(\phi^2)(t)| &\leq& H_{3\infty}L_4L_k t_0(t-t_0)|\phi^1-\phi^2|_{t_0t},
\\
|F_4(\phi^1)(t)-F_4(\phi^2)(t)| &\leq & L_2L_k(t-t_0)^2|\phi^1-\phi^2|_{t_0t},
\\
|F_5(\phi^1)(t)-F_5(\phi^2)(t)| &\leq & H_{3\infty}L_4L_k(t-t_0)^2|\phi^1-\phi^2|_{t_0t}+H_{4\infty}L_3(t-t_0)|\phi^1-\phi^2|_{t_0t}.
\end{eqnarray*}
Hence combining the above expressions we get, for $t\in[t_0, T]$
$$|F(\phi^1)(t)-F(\phi^2)(t)|\leq (L_1L_k+L_2L_k T+H_{3\infty}L_4L_k T+H_{4\infty}L_3 )(t-t_0)|\phi^1-\phi^2|_{t_0t}.$$
Hence we can choose $h>0$ small enough so that  $ \sup_{[t_0,t]}|F(\phi^1)(s)-F(\phi^2)(s)|<\varrho |\phi^1-\phi^2|_{t_0t}$
for some positive $\varrho<1$ and $t-t_0=h$. So by contraction mapping theorem, there exists a unique continuous function $x$ defined on $[t_0, t]$
satisfying \eqref{apen-a2}.

Putting $t_0=0$, we see that $x_t$ satisfies \eqref{appn-a1} on $[0, h]$. Having the solution defined on $[0, nh\wedge T]$, we can extend it uniquely on
$[0, (n+1)h\wedge T]$ for $n\in \N$. Since $h>0$ is fixed, this defines the solution uniquely on $[0,T]$.\hfill $\Box$.
\spa

\noi We can extend Lemma \ref{appn-lem1} as follows:

\begin{lemma}\label{appn-lem2}
Fix $T>0$. Let $k,H_i:\R\to\R, i=1,2,3,4,5,$ be Lipschitz continuous. We also assume $H_4, H_5$ to be bounded. Then the following integral equation
\begin{eqnarray}
x_t &=& H_1(S(0,t))+\int_0^tH_2(S(s,t))ds
 +\int_0^tH_3(x_s)H_5(x_s)H_4(S(s,t))ds\label{appn-a3}
\\
x_0 &=& H_1(0).\nonumber
\end{eqnarray}
has a unique solution in $C([0,T], \R)$ where $S(s,t)=\int_s^tk(x_u)du$.
\end{lemma}

\noi{\bf Proof:}\ Let $\varphi_n:\R\to\R$ be a smooth cut-off function such that $0\leq \varphi_n\leq 1,\ \varphi_n(s)=1$ on $[-n,n]$ and $\varphi_n(s)=0$ outside of $[-n-1,n+1]$.
Define $H^n_i(s)=\varphi_n(s)H_i(s)$ for $i=1,2,3$. Then $H^n_i$ is a bounded, Lipschitz continuous function for $i=1,2,3$. Hence $H^n_3H_5$ is a bounded, Lipschitz continuous
function for $n\in\N$. Therefore applying Lemma \ref{appn-lem1}, we have unique $x^n:[0,T]\to\R$, continuous, satisfying
\begin{eqnarray}
x^n_t &=& H^n_1(\int_0^tk(x^n_s)ds)+\int_0^tH^n_2(\int_s^tk(x^n_u)du)ds\nonumber
\\
&& +\int_0^tH^n_3(x^n_s)H_5(x^n_s)H_4(\int_s^tk(x^n_u)du)ds\label{appn-a4}
\\
x^n_0 &=& x_0=H_1(0).\nonumber
\end{eqnarray}
Since $k, H_i, i=1,\ldots, 5$ are Lipschitz and $H_4,H_5$ are bounded, we can get positive constants $d_1, d_2$ such that
\begin{eqnarray*}
 |H^n_1(\int_0^tk(x^n_s)ds)|&\leq & d_1+d_2\int_0^t|x^n_s|ds
 \\
 |\int_0^tH^n_3(x^n_s)H_5(x^n_s)H_4(\int_s^tk(x^n_u)du)ds| &\leq & d_1+ d_2\int_0^t|x^n_s|ds ,
\end{eqnarray*}
and
\begin{eqnarray*}
|\int_0^tH^n_2(\int_s^tk(x^n_u)du)ds| \leq  d_1+d_2 \int_0^t\int_s^t |x^n_u|du ds\leq d_1+d_2T\int_0^t |x^n_s| ds,
\end{eqnarray*}
for all $t\in[0, T]$ and $n\in\N$. Now combining these estimates with \eqref{appn-a4} and applying Gronwall's inequality we have
\be\label{appn-a5}
\sup_n\sup_{[0,T]}|x^n_s|\leq d_3,
\ee
for some constant $d_3$. For any compact $C\subset\R,$ there exists constant $L_C$ such that $|H^n_i(x)-H^n_i(y)|\leq L_c|x-y|$ for all $x,y\in C$ and $n\in \N$.
Therefore using an analogous expression as \eqref{apen-a2}, we get a constant $d_4>0$ (depending on $d_3$) satisfying
\be\label{appn-a6}
\sup_n\sup_{0\leq s\leq t\leq T}|x^n_s-x^n_t|\leq d_4|t-s|.
\ee
\eqref{appn-a5} and \eqref{appn-a6} imply that the sequence $\{x^n\}$ is equi-continuous family of continuous functions on $[0, T]$. Therefore using Arzel\'a-Ascoli, theorem there
is a $x:[0,T]\to\R$, continuous, such that $|x^{n_k}- x|_{0T}\to 0$ along some subsequence $n_k\to\infty$. Hence letting $n_k\to\infty$ in \eqref{appn-a4} we have
\begin{eqnarray*}
x_t &=& H_1(\int_0^tk(x_s)ds)+\int_0^tH_2(\int_s^tk(x_u)du)ds
\\
&&\ +\int_0^tH_3(x_s)H_5(x_s)H_4(\int_s^tk(x_u)du)ds
\\
x_0 &=& x_0=H_1(0).
\end{eqnarray*}
This proves the existence of solution to \eqref{appn-a3}. To prove uniqueness, let $\bar x$ be another solution to \eqref{appn-a3}. Define
$\sigma_n=\inf\{t\geq 0: |\bar x_t|>n\}\wedge T$. Since $H^n(s)=H(s)$ for $|s|\leq n$, from Lemma \ref{appn-lem1}, we get $x_s=x^n_s=\bar x_s$ for $s\leq\sigma_n$ and for all $n$
large (we need to take large $n$ to ensure that $x_0\in [-n,n]$). Therefore to complete the proof it is enough to show that $\liminf_{n\to\infty}\sigma_n=T$. But this is
obvious as $\sup_{[0,T]}|\bar x_s|<\infty$ (follows from a simple calculation similar to \eqref{appn-a5}).\hfill $\Box$

\subsection{Appendix B}
Consider a sequence of probability measures $\{\nu^n\}$ and $\nu$ on $[0, \infty)$ such that $\nu^n\to\nu$ as $n\to\infty$.
Let $\{v^n_i\}_{i=-\infty}^{i=\infty}$ be an i.i.d. sequence with common probability distribution $\nu^n$.
For $m\in\mathbb{Z}$ and $\ell\geq 0$, define
\be\label{a1}
\Li^n(m, \ell)=\frac{1}{n}\sum_{i=1+m}^{m+\lfloor n\ell\rfloor}\del_{v^n_i}.\tag{B1}
\ee
By Skorohod representation theorem, there exists $[0, \infty)$-valued random variables $Y^n\sim\nu^n$ and $Y\sim\nu$ such that $Y^n\to Y$ almost surely on some common probability space.
Define $\bar Y=\sup_n Y^n$. Let $\bar\nu$ be the law of $\bar Y$. There exists a continuous, increasing, unbounded function
$\bar f$ such that $\bar f\geq 1$ and $\lan\bar f^2, \bar\nu\ran<\infty$ (see Appendix B in \cite{zhang-dai-zwart}).
Define
$$\mathcal{V}=\{\chi_{C_x}, x\geq 0\}\cup\{\chi_{\bar C_x}, x\geq 0\}\cup\{\bar f\}.$$

\begin{lemma}\label{apen}
Fix $M, L>0$. Under the above assumptions, for all $\eps, \eta>0$ we have
\begin{displaymath}
\limsup_{n\to\infty} \p_n(\max_{-nM<m<nM}\sup_{\ell\in [0, L]}\sup_{f\in\mathcal{V}}|\lan f, \Li^n(m, \ell)\ran-\ell\lan f, \nu^n\ran|>\eps)<\eta.
\end{displaymath}
\end{lemma}
For the proof of above lemma we refer Lemma B.1 in \cite{zhang-dai-zwart}. Following the same argument as in Lemma 5.1 in \cite{zhang-dai-zwart}
we can have a sequence $\eps_A(n)$ such that $\eps_A(n)\to 0$ as $n\to\infty$ and $\lim_{n\to\infty}\p_n(\Omega_A^n(M, L))=1$ for every
fixed $M, L>0$ where

\be\label{a2}
\Omega_A^n(M, L)=\{\max_{-nM<m<nM}\sup_{\ell\in [0, L]}\sup_{f\in\mathcal{V}}|\lan f, \Li^n(m, \ell)\ran-\ell\lan f, \nu^n\ran|\leq\eps_A(n)\}.\tag{B2}
\ee

\spa

\noi{\bf Acknowledgement:} The author is grateful to Prof. Rami Atar for his valuable suggestions to improve this paper.

\bibliographystyle{plain}

\bibliography{refs}

\begin{thebibliography}{10}

\bibitem{atar-kaspi-shim}
Rami Atar, Haya Kaspi, and Nahum Shimkin.
\newblock Fluid limits for many server systems with reneging under a priority
  policy.
\newblock {\em Submitted}.

\bibitem{brown}
Lawrence Brown, Noah Gans, Avishai Mandelbaum, Anat Sakov, Haipeng Shen, Sergey
  Zeltyn, and Linda Zhao.
\newblock Statistical analysis of a telephone call center: a queueing-science
  perspective.
\newblock {\em J. Amer. Statist. Assoc.}, 100(469):36--50, 2005.

\bibitem{delay-verejones}
D.~J. Daley and D.~Vere-Jones.
\newblock {\em An introduction to the theory of point processes. {V}ol. {I}}.
\newblock Probability and its Applications (New York). Springer-Verlag, New
  York, second edition, 2003.
\newblock Elementary theory and methods.

\bibitem{dec-pascal}
L.~Decreusefond and P.~Moyal.
\newblock Fluid limit of a heavily loaded {EDF} queue with impatient customers.
\newblock {\em Markov Process. Related Fields}, 14(1):131--158, 2008.

\bibitem{ethier-kurtz}
Stewart~N. Ethier and Thomas~G. Kurtz.
\newblock {\em Markov processes}.
\newblock Wiley Series in Probability and Mathematical Statistics: Probability
  and Mathematical Statistics. John Wiley \& Sons Inc., New York, 1986.
\newblock Characterization and convergence.

\bibitem{gro-rob-zwart}
H.~Christian Gromoll, Philippe Robert, and Bert Zwart.
\newblock Fluid limits for processor-sharing queues with impatience.
\newblock {\em Math. Oper. Res.}, 33(2):375--402, 2008.

\bibitem{halfin-whitt}
Shlomo Halfin and Ward Whitt.
\newblock Heavy-traffic limits for queues with many exponential servers.
\newblock {\em Oper. Res.}, 29(3):567--588, 1981.

\bibitem{kang-ramanan}
Weining Kang and Kavita Ramanan.
\newblock Fluid limits of many-server queues with reneging.
\newblock {\em Ann. Appl. Probab.}, 20(6):2204--2260, 2010.

\bibitem{kang-ramanan1}
Weining Kang and Kavita Ramanan.
\newblock Asymptotic approximation for stationary distribution of many-server
  queues with abandonment.
\newblock {\em Ann. Appl. Probab.}, 22(6):477--521, 2012.

\bibitem{kaspi-ramanan}
Haya Kaspi and Kavita Ramanan.
\newblock Law of large numbers limits for many-server queues.
\newblock {\em Ann. Appl. Probab.}, 21(1):33--114, 2011.

\bibitem{liu-whitt}
Yunan Liu and Ward Whitt.
\newblock Large-time asymptotics for the {$G_t/M_t/s_t+GI_t$} many-server fluid
  queue with abandonment.
\newblock {\em Queueing Syst.}, 67(2):145--182, 2011.

\bibitem{mandel-mass-reim}
Avi Mandelbaum, William~A. Massey, and Martin~I. Reiman.
\newblock Strong approximations for {M}arkovian service networks.
\newblock {\em Queueing Systems Theory Appl.}, 30(1-2):149--201, 1998.

\bibitem{ramanan-reiman}
Kavita Ramanan and Martin~I. Reiman.
\newblock Fluid and heavy traffic diffusion limits for a generalized processor
  sharing model.
\newblock {\em Ann. Appl. Probab.}, 13(1):100--139, 2003.

\bibitem{reed}
Josh Reed.
\newblock The {$G/GI/N$} queue in the {H}alfin-{W}hitt regime.
\newblock {\em Ann. Appl. Probab.}, 19(6):2211--2269, 2009.

\bibitem{whitt}
Ward Whitt.
\newblock Fluid models for multiserver queues with abandonments.
\newblock {\em Oper. Res.}, 54(1):37--54, 2006.

\bibitem{yamada}
Keigo Yamada.
\newblock Diffusion approximation for open state-dependent queueing networks in
  the heavy traffic situation.
\newblock {\em Ann. Appl. Probab.}, 5(4):958--982, 1995.

\bibitem{zhang}
Jiheng Zhang.
\newblock Fluid models of many-server queues with abandonment.
\newblock {\em Queueing Systems Theory Appl.}, Forthcoming, 2012.

\bibitem{zhang-dai-zwart}
Jiheng Zhang, J.~G. Dai, and Bert Zwart.
\newblock Law of large number limits of limited processor-sharing queues.
\newblock {\em Math. Oper. Res.}, 34(4):937--970, 2009.

\end{thebibliography}

\end{document}